\documentclass[12pt,a4paper,leqno]{article}

\usepackage[utf8]{inputenc}
\usepackage{amsmath}
\usepackage{amsfonts}
\usepackage{amsthm}
\usepackage{calrsfs}
\usepackage{amssymb}

\usepackage{dsfont}
\usepackage{mathtools}
\usepackage{color}
\usepackage{float}

\usepackage{geometry}
\newtheorem{theorem}{Théorème}%%compteur

\newtheorem{prop}[theorem]{\textbf{Proposition}}

\newenvironment{dem}{\textbf{Proof}}{\qed \\ }
\newenvironment{demof}{\textbf{Proof of }}{\qed \\ }

\newtheorem{defi}[theorem]{\textbf{Definition}}
\newtheorem{lemm}[theorem]{\textbf{Lemma}}
\newtheorem{thm}[theorem]{\textbf{Theorem}}

\newtheorem*{thm*}{\textbf{Theorem}}%thm sans numero
\newtheorem{assump}{\textbf{Assumption}}

\usepackage{filecontents,lipsum}

\usepackage{mystyle}

\newcommand{\R}{\mathbb{R}}

\newcommand{\E}{\mathbb{E}}

\newcommand{\cL}{\mathcal{L}}
\title{Behavior of the Poincaré constant along the Polchinski renormalization flow}
\author{Jordan Serres\footnote{Institut de Mathématiques de Toulouse, jordan.serres@math.univ-toulouse.fr}}

\begin{document}
\maketitle

\begin{center}
\textbf{Abstract}
\end{center}
\begin{quotation}
We control the behavior of the Poincaré constant along the Polchinski renormalization flow using a dynamic version of $\Gamma$-calculus. We also treat the case of higher order eigenvalues. Our method generalizes a method introduced by B. Klartag and E. Putterman to analyze the evolution of log-concave distributions along the heat flow. Furthermore, we apply it to general $\varphi^4$-measures and discuss the interpretation in terms of transport maps.
\end{quotation}

\section{Introduction}

A probability measure $\mu$ on $\mathbb{R}^d$ is said to satisfy a Poincaré inequality when there exists a positive finite constant $C$ such that for all functions $f$ in the Sobolev space $H^{1}(\mu)$,
$$
\mathrm{Var}_\mu (f) \leq C \int |\nabla f|^2 d\mu.  
$$
We denote by $C_P(\mu)$ the smallest constant for which the above inequality holds.
Poincaré inequalities have many applications (see for instance \cite{surveyineg}), and are in particular related to the Kannan-Lovasz-Simonovits conjecture which states that all isotropic log-concave distributions have dimension-independent Poincaré constants (see \cite{surveykls} for a survey).\\
A classical method to obtain bounds on the Poincaré constant $C_P(\mu)$ of a distribution $\mu$ is the perturbative approach: if $C_P(\mu)$ is known, we want to compare it with the Poincaré constant $C_P(\delta\mu)$ of a small perturbation $\delta\mu$. %This problem is the converse of the problem of stability which want to bound the distance between $\mu$ and $\nu$ when their Poincaré constants $C_P(\mu)$ and $C_P(\nu)$ are closed.
Some results exist when $\delta\mu$ is a bounded perturbation of $\mu$, i.e. when $\delta\mu$ has a bounded density with respect to $\mu$ (see \cite[Proposition 4.2.7]{BGL}). Another result is that if $\gamma$ denotes the standard Normal distribution on $\R^d$ (then $C_P(\gamma)=1$) and $\ast$ denotes the convolution, then for all probability distributions $\mu$ we have $C_P(\mu)\geq C_P(\mu\ast\gamma)-C_P(\gamma)$. This is a consequence of \cite[Proposition 1]{ballbarthenaor}. The reverse inequality 
\begin{equation}\label{cattiauxguillin}
C_P(\mu)\leq C_P(\mu\ast\gamma)+C_P(\gamma) 
\end{equation}
is more informative and was proved under a log-concavity assumption on the measure $\mu$ (see \cite[Theorem 9.4.3]{cattiauxguillin}). In \cite{KP}, B. Klartag and E. Putterman interpret the convolution with the Gaussian distribution as the action of the heat flow on $\mu$. This interpretation together with the preservation of the log-concavity by the heat flow, allowed them to strengthen $\eqref{cattiauxguillin}$ and to prove by $\Gamma$-calculus that the Poincaré constant is non-decreasing along the heat flow acting on a log-concave distribution.
\begin{thm*}\cite[Theorem 1.1]{KP}
If $\mu$ is log-concave, then for all $0\leq s\leq t$, 
\begin{equation}\label{thmkp}
C_P(\mu\ast\gamma_s)\leq C_P(\mu\ast\gamma_t), 
\end{equation} where $\gamma_s$ denotes the Normal distribution with covariance matrix $sI_d$.
\end{thm*}
They also proved $\eqref{thmkp}$ by constructing contractive transport maps, see Section $\ref{sectiontransport}$ for a discussion of this method.
Let us mention that, following this, the analysis of the evolution of a log-concave distribution along the heat flow allowed B. Klartag and J. Lehec to obtain the best result currently known towards the KLS conjecture (see \cite{klartaglehec}).\\
In this paper, we study the behavior of the Poincaré constant along the Polchinski flow. The origins of this flow go back to the work of J. Polchinski \cite{polchinski1984renormalization} on the Renormalization Group method. It was introduced in its current formulation with the multiscale Bakry-Emery criterion by R. Bauerschmidt and T. Bodineau in \cite{BB} in order to derive logarithmic-Sobolev inequality for the continuum sine-Gordon model. Subsequently, this flow has also been used to derive logarithmic-Sobolev inequalities for other models from Quantum Fields Theory, such as the $\varphi^4$ models \cite{bauerschmidt2022log} or the Ising models \cite{bauerschmidt2022ising}.\\ 
Let $\nu_0= e^{-V_0(x)}\gamma_{C_\infty}$ be a probability distribution with density $e^{-V_0}$ with respect to the Gaussian distribution $\gamma_{C_\infty}$ with covariance (positive definite) matrix $C_\infty$. A family $(C_t)_{t\geq 0}$ of semidefinite matrices on $\R^d$ increasing continuously as quadratic forms from $C_0=0$ to $C_\infty$ is called a covariance decomposition. Given a covariance decomposition, we can define the renormalized potential
$$V_t  := -\log \left( \gamma_{C_t} \ast e^{-V_0}\right).$$ The Polchinski flow is then defined as $$\nu_t  := e^{V_\infty(0)-V_t(x)}\gamma_{C_\infty-C_t}.$$
By definition, $\gamma_0:=\delta_0$ is the Dirac delta function, hence $\nu_\infty=\delta_0$ and the Polchinski flow interpolates between $\nu_0$ and $\delta_0$.
We prove the following theorem.
\begin{thm*}
If the Polchinski flow $(\nu_t)_t$ of a probability distribution $\nu_0$ satisfies the multiscale Bakry-Emery criterion, then the Poincaré constant along the flow satisfies
$$
\forall 0<s\leq t,\quad C_P(\nu_s)\leq e^{(\alpha_t-\alpha_s)-2(\lambda_t - \lambda_s)}C_P(\nu_t), 
$$
where $\lambda_t$ are the multiscale curvatures, and $\alpha_t$ are corrective terms due to the fact that the Polchinski flow admits the Dirac delta function as ergodic distribution (see Definition $\ref{defalpha}$).
\end{thm*}

This paper is organized as follows. In Section $\ref{sectiongammacalc}$, we introduce the general $\Gamma$-calculus method to study the behavior of the Poincaré constant along a flow of probability distributions driven by a diffusion semigroup. We also treat the case of higher order eigenvalues. In Section $\ref{sectiontransport}$, we discuss the transport map method, and in particular why it cannot currently be applied to the Polchinski flow. In Section $\ref{mainsection}$, we rigorously define the Polchinski flow, the multiscale Bakry-Emery criterion, the renormalized Poincaré constant along the flow, and then present a proof of the main theorem of this paper. Finally, in Section $\ref{sectionphi4}$, we illustrate our theorem on the example of general $\varphi^4$ measures.

\section{A general $\Gamma$-calculus approach}\label{sectiongammacalc}

Let $\cL$ be the generator of a Markov diffusion semigroup $P_t$ with carré du champ operator $\Gamma$ and reversible measure $\mu$ on a manifold $M$. Let $\mu_0$ be a probability distribution on $M$ with density $\rho$ with respect to $\mu$. If we draw randomly $X_0$ with distribution $\mu_0$ and let the semigroup evolve, then we obtain a Markov process $X_t$, whose laws $\mu_t$ interpolate between $\mu_0$ and $\mu_\infty=\mu$. In this section, we describe a general method for obtaining bounds on the Poincaré constant along the flow $\mu_t$ of the form
\begin{equation}\label{cpalongtheflowgen1}
\forall 0\leq s\leq t,\quad C_P(\mu_s)\leq e^{\int_s^t k_s'\,ds}C_P(\mu_t), 
\end{equation}
where $k_t'$ are curvature terms that will be described later. This monotonicity-type result for the Poincaré constant can be captured from the point of view of $\Gamma$-calculus. When the operator $\cL$ has a discrete spectrum, this method also allows to treat the case of higher order eigenvalues. The measure $\mu_t$ has density $P_t(\rho)$ with respect to $\mu$, and allows us to define the dual semigroup 
\begin{equation}\label{dualsemigrp}
Q_t(\phi):=\frac{P_t(\phi\rho)}{P_t(\rho)}.
\end{equation} This dual semigroup evolves in such a way that there exist a flow of diffusion operators $(L_t)_t$ satisfying $$\partial_t Q_t = L_tQ_t.$$ 
We can see from a straightforward computation that 
\begin{equation}\label{ltformula}
L_tu = \cL u + 2\,\Gamma\left(\log P_t(\rho),u\right).
\end{equation}
We denote by $\Gamma_t(f):=\frac{1}{2}(L_t(f^2)-2fL_tf)$ and $\Gamma_{2,t}:=\frac{1}{2}\left(L_t\Gamma_t(f) -2\Gamma_t(f,L_tf) \right)$ its iterated first and second $\Gamma$-operators. It is immediate to see that $\Gamma_t=\Gamma$ does not depend on time. We also define the operators $\Lambda_t$ as the diffusion generators with reversible distribution $\mu_t$, and denote by $\Gamma_t^\Lambda(f):=\frac{1}{2}(\Lambda(f^2)-2f\Lambda_tf)$ and $\Gamma_{2,t}^\Lambda=\frac{1}{2}\left(\Lambda_t\Gamma_t^\Lambda(f) -2\Gamma_t^\Lambda(f,\Lambda_tf) \right)$ its iterated first and second $\Gamma$-operators. An easy computation shows that
\begin{equation}\label{lambdatop}
\Lambda_tu = 2\cL u + 2\,\Gamma\left(\log P_t(\rho),u\right),
\end{equation}
and moreover, $\Gamma_t^\Lambda=2\,\Gamma$ does not depend on time. Although in this case none of the operators $\Gamma_t$ and $\Gamma_t^\Lambda$ are time dependent, we still retain the notation in order to better understand the method we will use for the Polchinski flow where the operators will indeed be time dependent.
The key strategy to obtain $\eqref{cpalongtheflowgen1}$ is then to show the quasi-decay of the Raylegh quotient $$R_\phi(t) := \frac{\E_{\mu_t} \left[\Gamma_t^\Lambda\left( Q_t\phi\right)\right]}{\E_{\mu_t}\left[(Q_t\phi)^2\right]}, $$ when $\phi : M\rightarrow\R$ is a smooth and compactly supported function. For this purpose, using that $\Gamma_t^\Lambda=2\,\Gamma_t $, we can compute that
$$
\partial_t \E_{\mu_t}\left[(Q_t\phi)^2\right] = -2\,\E_{\mu_t}\left[\Gamma_t(Q_t\phi)\right]= -\E_{\mu_t} \left[\Gamma_t^\Lambda\left( Q_t\phi\right)\right]\quad\mathrm{and}\quad \partial_t \E_{\mu_t} \left[\Gamma_t^\Lambda\left( Q_t\phi\right)\right] = 4\,\E_{\nu_t}\left[ \Gamma_{2,t}(Q_t\phi)\right],
$$
and therefore
\begin{align*}
R_\phi'(t) &= \frac{\left( E_{\mu_t}\left[\frac{d\Gamma_t^\Lambda}{dt}(Q_t\phi)\right]-4\,\E_{\mu_t}\left[ \Gamma_{2,t}(Q_t\phi) \right]\right)\E_{\mu_t}\left[(Q_t\phi)^2\right]+ \E_{\mu_t} \left[\Gamma_t^\Lambda\left( Q_t\phi\right)\right]^2 }{\E_{\mu_t}\left[(Q_t\phi)^2\right]^2}\\
&= \frac{ \E_{\mu_t} \left[\Gamma_t^\Lambda\left( Q_t\phi\right)\right]^2 -\E_{\mu_t}\left[ \Gamma_{2,t}^\Lambda(Q_t\phi) \right]\E_{\mu_t}\left[(Q_t\phi)^2\right] }{\E_{\mu_t}\left[(Q_t\phi)^2\right]^2} + \frac{E_{\mu_t}\left[\frac{d\Gamma_t^\Lambda}{dt}(Q_t\phi)\right] + \E_{\mu_t}\left[ \left(\Gamma_{2,t}^\Lambda - 4\,\Gamma_{2,t}\right)(Q_t\phi) \right]}{\E_{\mu_t}\left[(Q_t\phi)^2\right]}.
\end{align*}
We can easily see that the term $\E_{\mu_t} \left[\Gamma_t^\Lambda\left( Q_t\phi\right)\right]^2 -\E_{\mu_t}\left[ \Gamma_{2,t}^\Lambda(Q_t\phi) \right]\E_{\mu_t}\left[(Q_t\phi)^2\right]$ is non positive. Indeed, we have by integrating by parts and using the Cauchy-Schwarz inequality
\begin{align*}
\E_{\mu_t} \left[\Gamma_t^\Lambda\left( Q_t\phi\right)\right] & = -\E_{\mu_t}\left[ Q_t\phi (\Lambda_t Q_t\phi) \right] \\
& \leq \sqrt{\E_{\mu_t}\left[ (Q_t\phi)^2\right] \E_{\mu_t}\left[ (\Lambda_t Q_t\phi)^2 \right]}\\
& = \sqrt{\E_{\mu_t}\left[(Q_t\phi)^2\right] \E_{\mu_t}\left[ \Gamma_{2,t}^\Lambda(Q_t\phi) \right]}.
\end{align*}
It appears therefore that under the dynamical curvature assumption
\begin{equation}\label{cddynamic}
\frac{d\Gamma_t^\Lambda}{dt} + \left(\Gamma_{2,t}^\Lambda - 4\,\Gamma_{2,t}\right)  \leq  k_t' \Gamma_t,
\end{equation}
the Gronwall lemma gives us that the Raylegh quotient satisfies $R_\phi(t)\leq e^{\int_0^t k_s'ds}R_\phi(0)$ and this is enough to show $\eqref{cpalongtheflowgen1}$.\\
When the operator $\cL$ has discrete spectrum, let us denote by $\lambda_k(\mu)$, $k\geq 1$ its $k$-th non-zero eigenvalue. By the min-max theorem, there exists a ($k+1$)-dimensional subspace $E\subset H^1(\mu)$ such that $\forall \phi \in E\setminus \{0\}$, $R_\phi(0)\leq \lambda_k(\mu)$. If furthermore, for all $t>0$, the dual operator $Q_t$ given by $\eqref{dualsemigrp}$ is a bijection from $L^1(\mu)$ onto $L^1(\mu)$, then the space $E_t:=\{Q_t(\phi)\,|\,\phi\in E\}$ is a ($k+1$)-dimensional vector space in $H^1(\mu_t)$ such that for all $\psi=Q_t\phi\in E_t$, $$\frac{\E_{\mu_t} \left[\Gamma_t^\Lambda\left( \psi\right)\right]}{\E_{\mu_t}\left[\psi^2\right]} = R_\phi(t) \leq e^{\int_0^t k_s'ds}R_\phi(0)\leq e^{\int_0^t k_s'ds} \lambda_k(\mu).$$ Again by the min-max theorem, we can conclude that $$\lambda_k(\mu_t)\leq e^{\int_0^t k_s'ds}\lambda_k(\mu).$$

When the operator $\cL$ is of the form $\Delta - \nabla \Phi\cdot \nabla$, then $L_tu= \Delta u + \nabla\left(2\log P_t\rho -\Phi\right)\cdot\nabla u $, $\Lambda_t u = 2\,\Delta u + 2\,\nabla\left(\log P_t\rho -\Phi\right)\cdot\nabla u $, $\Gamma_t(u)=|\nabla u|^2$, $\Gamma_t^\Lambda(u)=2\,|\nabla u|^2$, and $\Gamma_{2,t}^\Lambda = 4\,\Gamma_{2,t} + \nabla^2 \log P_t(\rho)$. Consequently, the condition $\eqref{cddynamic}$ becomes the classical Bakry-Emery condition $\nabla^2 \log P_t(\rho) \leq k' I_d $.\\

In \cite{KP}, B. Klartag and E. Putterman take the operator $\cL$ as (half of) the Laplacian on $\R^d$. Therefore the reversible distribution $\mu$ is the Lebesgue measure on $\R^d$ and the semigroup $P_t$ is the heat semigroup. They start from a log-concave distribution $\mu_0$ with Lebesgue density $\rho$, and by stopping the flow at time $t=1$, they obtain an interpolation between $\mu_0$ and $\mu_1=P_1(\rho)dx=\mu_0\ast \gamma$ where $\gamma$ denotes the normalized Gaussian distribution on $\R^d$ and $\ast$ denotes the convolution. The dynamical curvature assumption $\nabla^2 \log P_t(\rho)\leq 0 $ is then satisfied because the heat flow preserves the log-concavity. Hence $k_t'$ can be chosen to be zero and $\eqref{cpalongtheflowgen1}$ becomes an exact monotonicity result. Moreover, the dual operator $Q_t$ is indeed a bijection, since the quantity $P_t(\phi\rho)= \phi\rho\ast\gamma_t$ allows to recover the Fourier transform of $\phi\rho$. Therefore the monotonicity result is satisfied by all eigenvalues of any order.\\

If we take the Ornstein-Ulhenbeck generator $\cL = \Delta - A\varphi\cdot \nabla$, then the reversible distribution is the Normal distribution $\gamma_{A^{-1}}$ with covariance matrix $A^{-1}$, and the semigroup is given by $P_t(f)(\varphi) = \E_{\gamma_{C_{2t}}}\left[f\left(e^{t A}\varphi +\zeta \right) \right] $ where $C_t= A^{-1}-A^{-1}e^{-tA}$. Consider then a probability distribution $\mu_0$ which is absolutely continuous with respect to $\gamma_{A^{-1}}$ and with density $e^{-V_0}$. Then $(\mu_t)_t$ is the Langevin flow which interpolates between $\mu_0$ and the Gaussian $\gamma_{A^{-1}}$, and the density of $\mu_t$ with respect to $\gamma_{A^{-1}}$ is given by $P_\frac{t}{2}\left(e^{-V_0}\right)(\varphi) = e^{-V_t(e^{\frac{t}{2}A}\varphi}) $, where $V_t$ is the renormalized potential which we will define in Section $\ref{defrenorm}$ (see \cite[Lemma 2.1]{shenfeld2022exact}). 
Therefore $\nabla^2\log P_{\frac{t}{2}}\left(e^{-V_0}\right) = e^{\frac{t}{2}A}\nabla^2 V_t e^{\frac{t}{2}A} $, and the dynamical curvature assumption $\nabla^2 \log P_t(\rho) \leq k_t' I_d $ is no other than the multiscale Bakry-Emery criterion $\eqref{multibe2}$. We have then showed that if $\mu_0=e^{-V_0}d\gamma_{A^{-1}}$ satisfies the multiscale Bakry-Emery criterion $\eqref{multibe2}$, then the Poincaré constant along the Langevin flow satisfies $\eqref{cpalongtheflowgen1}$. Moreover, the dual operator $Q_t(\phi)=\frac{P_t(e^{-V_0}\phi)}{P_t (e^{-V_0})}$ is again a bijection on $L^1(\gamma_{A^{-1}})$. Therefore the quasi monotonicity result is satisfied by all eigenvalues of any order. These results can be also proved from \cite[Theorem 2.2]{shenfeld2022exact}, where Y. Shenfeld shows that under the multiscale Bakry-Emery criterion $\eqref{multibe2}$, the Langevin transport map constructed as we will see in Section $\ref{sectiontransport}$ is Lipschitz.\\

In the case of the Polchinski flow $(\nu_t)_t$ (see Section $\ref{defrenorm}$), the previous framework does not apply. Indeed, this framework would apply to the flow $(\mu_t)_t$ where $d\mu_t(x)=e^{V_\infty(x)-V_t(x)}dx$, however the distributions $\nu_t$ have in addition the Gaussian term $e^{-\frac{1}{2}<x,\,(C_\infty - C_t)^{-1}x>}$. Therefore, in order to study the behavior of the Poincaré constant along the flow, we have to introduce in addition to $L_t$ and $\Lambda_t$, the diffusion generators $\cL_t$ defined in $\eqref{generateurdenut}$. As a result, we can derive a monotonicity-type formula such as $\eqref{cpalongtheflowgen1}$, but only at the cost of an additional term $e^{(\alpha_t-\alpha_s)}$ depending on the Gaussian factor $e^{-\frac{1}{2}<x,\,(C_\infty - C_t)^{-1}x>}$ (see Theorem $\ref{thmmonotonicity}$).

\section{A transport map interpretation}\label{sectiontransport}

When a probability measure $\nu$ is the pushforward by a Lipschitz function $S$ of a probability distribution $\mu$ which satisfies some functional inequalities (such as the Poincaré inequality or the logarithmic-Sobolev inequality), then these functional inequalities are transfered to $\nu$ via the map $S$ (see \cite{cordero-erosquinpushforward, Milman2018}). Consequently, a powerful way to derive functional inequalities for a probability distribution is to construct a transport map from a classic probability distribution satisfying those functional inequalities. These are generally log-concave distributions (see \cite{kimmilman}), the Normal distribution (see e.g \cite{transportviaheatflow,neeman2022lipschitz}), or the Wiener measure (see \cite{browniantransportmap}). Put in a general and informal way, the conditions on the measure $\nu$ to be realized as a transport of one of the previous measures, are a compromise between the smallness of its support and the magnitude of its curvature (as the convexity of its potential). The Bakry-Emery multiscale criterion $\eqref{multibe2}$ has recently been used as a sufficient condition for a measure $\nu$ to be the pushforward by a Lipschitz function of a Normal distribution, and this has been applied in particular to the case of the sine-Gordon model (see \cite{shenfeld2022exact}). Let us briefly describe a classic method to construct such a Lipschitz transport map.

Let $\cL$ be a generator of a Markov diffusion semigroup $P_t$ with reversible measure $\mu$ on a manifold $M$. If we draw randomly $X_0$ with distribution $\mu_0$ and let the semigroup evolve, then we obtain a Markov process $X_t$, whose laws $\mu_t$ interpolate between $\mu_0$ and $\mu_\infty=\mu$. If we find the vector field $(\xi_t)_t$ driving the flow $\mu_t$, i.e. such that the continuity equation $$\partial_t\mu_t +\nabla\cdot\left(\xi_t\,\mu_t\right)=0 $$ holds (where $\nabla\cdot$ denotes the divergence operator), then we can compute the integral curves $(T_t:M\rightarrow M)_t$ of the velocity field by solving $$\partial_t T_t(x) = \xi(t, T_t(x)). $$ If all goes well (see \cite[Chapter 4]{santabrogliooptimalbook} for rigorous details), then for all $t\geq 0$, $\mu_t$ is the pushforward of $\mu_0$ by $T_t$ and has density $P_t(\rho)d\mu$, where $\rho$ denotes the density of $\mu_0$ with respect to $\mu$. Moreover, for all $t\geq 0$, $T_t$ are bijective with inverse functions $S_t$, and therefore $S:=\underset{t\rightarrow+\infty}{\lim} S_t$ is a transport map from $\mu_\infty=\mu$ to $\mu_0$. It was shown that this transport map is different than the Brenier optimal transport map (see \cite{tanana}). The curvature type assumptions on $\mu_0$ are then used to show that $S$ is Lipschitz via a control of the driving vector field $\xi_t$.

If the transport maps $S_t : \mu\rightarrow\mu_t$ are $e^{\int_0^t k_s'\,ds} $-Lipschitz for some function $k':\R_+\rightarrow\R$, then the transport map $S:\mu\rightarrow\mu_0$ is $e^{\int_0^\infty k_s'\,ds}$-Lipschitz (see \cite[Lemma 1]{transportviaheatflow}). In particular, by taking the composition $ T_t\circ S_s$ for $s<t$, we obtain a $e^{\int_s^t k_s'\,ds} $-Lipschitz transport map from $\mu_s$ to $\mu_t$. In term of Poincaré constants, it implies the following control along the flow: 
\begin{equation}\label{cpalongtheflowgen}
C_P(\mu_s)\leq e^{\int_s^t k_s'\,ds}C_P(\mu_t). 
\end{equation}
We recognize the same result as $\eqref{cpalongtheflowgen1}$ in Section $\ref{sectiongammacalc}$. This transport map approach also allows to derive such results for higher order eigenvalues.\\ 

In the case of the Polchinski flow $(\nu_t)_t$ defined in Section $\ref{defrenorm}$, the distributions $\nu_t$ interpolate between $\nu_0$ given in $\eqref{nuzero}$ and $\nu_\infty=\delta_0$ the Dirac delta function. Therefore, this dynamic cannot give a transport map, since $\nu_0$ cannot be realized as the pushforward of the Dirac delta function. However, a natural question is to find a rescaling of the Polchinski semigroup such that the ergodic distribution $\nu_\infty$ is no longer the Dirac delta function, so a transport map between $\nu_0$ and $\nu_\infty$ could be constructed. The theorem $\ref{thmmonotonicity}$, linked to the monotonicity result $\eqref{cpalongtheflowgen}$ valid when the transport maps exist and are Lipschitz, strongly suggests that such a rescaling is possible. Moreover, this rescaling must remove the term $e^{(\alpha_t-\alpha_s)}$ in $\eqref{ineqthmprincipal}$, and thus the largest eigenvalue $\alpha_t'$ of the matrix $(C_t')^\frac{1}{2} \left(\nabla^2V_t + (C_\infty -C_t)^{-1}\right)(C_t')^\frac{1}{2}$ should characterize the order of magnitude of the rescaling. However, the problem of finding exactly such a rescaling is still open.

\section{The Poincaré constant along the Polchinski flow}\label{mainsection}

\subsection{The Polchinski renormalization semigroup}\label{defrenorm}

In this section, we introduce the Polchinski renormalization semigroup. Let $(C_t)_{t\geq 0}$ be a family of semidefinite matrices on $\R^d$. We assume that the $C_t$ matrices increase continuously as quadratic forms from $C_0=0$ to a matrix $C_\infty$. We also assume the family to be twice differentiable with respect to the parameter $t$, and moreover that its first derivative $t\mapsto C_t'$ is bounded. We denote by $\gamma_{C_t}$ the (possibly degenerate) Gaussian measure with covariance $C_t$ and mean zero. Let 
\begin{equation}\label{nuzero}
\nu_0= e^{-\frac{1}{2}<x,C_\infty^{-1} x>-V_0(x)}dx
\end{equation}
be a probability distribution on $\mathbb{R}^d$, where $V_0:\R^d\rightarrow\R$ is a smooth function. The part $e^{-<x,C_\infty^{-1} x>}$ represents the log-concave part of the measure, while $e^{-V_0}$ represents the non log-concave part of $\nu_0$. %(If $\nabla^2 V_0\geq 0$ then $\nu_0$ is log-concave).
We define
\begin{align*}
V_t & := -\log \left( \gamma_{C_t} \ast e^{-V_0}\right)\\
P_{s,t}f & := e^{V_t}\, \gamma_{C_t-C_s}\ast fe^{-V_s}\\
\nu_t & := e^{V_\infty(0)-V_t(x)}\gamma_{C_\infty-C_t}= e^{-\frac{1}{2}<x,\,(C_\infty - C_t)^{-1}x> - V_t(x)+V_\infty(0)}dx
\end{align*}
where $\ast$ denotes the convolution. By definition, $\gamma_0:=\delta_0$, hence $\nu_\infty=\delta_0$.
As in \cite[Section 2.1]{BB}, we assume that $\E_{\nu_t}g\left(P_{0,t}f\right)$ is continuous in $t$, so that $$\E_{\nu_t}g\left(P_{0,t}f\right)\underset{t\rightarrow\infty}{\rightarrow}g\left(\E_{\nu_0}(f)\right). $$
The renormalized potential $V_t$ satisfies the Polchinski equation (see \cite[Proposition 2.1]{BB})
$$\partial_t V_t = \frac{1}{2}\left( \frac{1}{2}\Delta_{C_t'} + L_t \right) V_t, $$
where $$L_t := \frac{1}{2}\Delta_{C_t'} - <\nabla V_t, \nabla >_{C_t'}, $$
and the index $C_t'$ in the Laplacian or the dot product denotes that these operations are computed with respect to $C_t'$ on $\R^d$, i.e.
$$<U, V >_{C_t'}\,:=\sum_{i,j} (C_t')_{i,j}U_i V_j\quad \mathrm{and}\quad \Delta_{C_t'} f := \sum_{i,j} (C_t')_{i,j}\frac{\partial^2 f}{\partial x_{i,j}}. $$
The diffusion operator $\Lambda_t:=\frac{1}{2}\Delta_{C_t'} + L_t$ admits $e^{-V_t(x)}dx$ as reversible distribution. Let us define the diffusion operator 
\begin{equation}\label{generateurdenut}
\cL_t:=\Delta_{C_t'} - <\nabla V_t + (C_\infty-C_t)^{-1}x, \nabla >_{C_t'},
\end{equation}
which admits $\nu_t$ as reversible measure. The two operators $\Lambda_t$ and $\cL_t$ have the same carré du champ operator $|\nabla \cdot |_{C_t'}^2$, but different $\Gamma_2$-operators. We will denote by $\Gamma_{2,t}$ the $\Gamma_2$-operator associated with $L_t$ and by $\Gamma_{2}^{\cL_t}$ the one associated with $\cL_t$. These operators will be used in the proof of Lemma $\ref{lemmdecrraylegh}$.

\subsection{The multiscale Bakry-Emery criterion}
The multiscale Bakri-Emery criterion has been introduced in \cite{BB} in order to derive logarithmic Sobolev inequalities for certain non-log-concave probability distributions from Quantum Field Theory on which the classical Bakry-Emery criterion is unusable. One of these examples (the $\varphi^4$ model) is treated in Section $\ref{sectionphi4}$. 
There are several different multiscale Bakry-Emery conditions, depending on the choice of the covariance decomposition $(C_t)_t$ in the Polchinski renormalization. We will consider the following one, which is the most general.
\begin{assump}(Multiscale Bakry-Emery criterion)\label{cdgenralized}
For all $t\geq 0$, there exists $\lambda_t'\in\R$ (possibly negative) such that
\begin{equation}\label{multiscaleeq}
\forall x\in\R^d,\quad C_t'\,\nabla^2 V_t\,C_t' \geq \frac{1}{2}	C_t'' + \lambda_t' C_t',
\end{equation} in the sense of quadratic forms. %The function $t\mapsto\lambda_t'$ is assumed to be locally integrable.
\end{assump}
The criterion $\eqref{multiscaleeq}$ is satisfied for example by the $\varphi^4$ models under the Pauli-Villar decomposition (see Section $\ref{sectionphi4}$).
When $V_0=0$, then $V_t=0$, and one takes $C_t=C_\infty-C_\infty e^{-tC_\infty^{-1}}$. Then $C_t''=-C_\infty^{-1} C_t'$ is negative semi definite, and then Assumption $\ref{cdgenralized}$ is satisfied with $\lambda_t' \leq \frac{\lambda}{2}$ where $\lambda$ is the smallest eigenvalue of $C_\infty$.
When $\nabla^2 V_0\geq 0$, then $\nabla^2 V_t \geq 0$ for all $t\geq 0$ (see \cite[Example 1.3]{BB}). Hence by taking again $C_t=C_\infty-C_\infty e^{-tC_\infty^{-1}}$ as in the case $V_0=0$, Assumption $\ref{cdgenralized}$ is satisfied again with $\lambda_t' \leq \frac{\lambda}{2}$.
The second and most classic multiscale Bakry-Emery criterion requires that 
\begin{equation}\label{multibe2}
 (C_t')^\frac{1}{2}\mathrm{Hess}V_t(C_t')^\frac{1}{2}\geq \mu_t'\,I_d 
\end{equation} and is more adapted to the heat kernel decomposition $C_t'=e^{-tC_\infty^{-1}}$, where in that case it implies $\eqref{multiscaleeq}$. This criterion is satisfied for example by the sine-Gordon model. 

\subsection{The Poincaré inequality under renormalization}

Under the multiscale Bakry-Emery criterion, the probability distribution $\nu_0$ satisfies a Poincaré inequality with sharp constant 
\begin{equation}\label{piundermbe}
C_P(\nu_0)\leq |C_0'|\int_0^\infty e^{-2\lambda_t}dt,
\end{equation} where $|C_0'|$ denotes the spectral radius of the matrix $C_0'$. Let us sketch the proof of this fact. First, we decompose the variance under $\nu_0$ along the renormalization flow $\nu_t$: 
\begin{equation}\label{vardecomp}
\mathrm{Var}_{\nu_0} (F) = -\int_0^\infty \frac{\partial}{\partial t}\E_{\nu_t}\left[(P_{0t}F)^2\right]\,dt =\int_0^\infty \E_{\nu_t}\left|\nabla P_{0t}F\right|_{C_t'}^2 dt.
\end{equation}
Next, it is possible to show that the multiscale Bakry-Emery criterion implies, as in the classical case, the following intertwining relation between the Polchinski semigroup and the gradient operator:
\begin{equation}\label{intertwining}
\left|\nabla P_{0t}F\right|_{C_t'}^2\leq |C_0'|e^{-2\lambda_t}P_{0t}\left(|\nabla F|^2\right).
\end{equation} 
Then $\eqref{vardecomp}$ and $\eqref{intertwining}$ together give $\eqref{piundermbe}$.
The proof of $\eqref{intertwining}$ is based on classic semigroup calculations. Let $t>0$ and let us define $\psi(s)=e^{-2(\lambda_t-\lambda_s)}P_{s,t}\left[\left|\nabla P_{0,s}F \right|_{C_s'}^2 \right]$ for all $s\in[0,t]$. One can compute that (see \cite[Lemma 2.8]{BB}) $$\psi'(s)= e^{-2(\lambda_t-\lambda_s)}P_{s,t}\left[2\lambda_s'\left|\nabla P_{0,s}F \right|_{C_s'}^2 - \left( L_s-\partial_s \right)\left(\left|\nabla P_{0,s}F \right|_{C_s'}^2\right)\right]\leq 0.$$ Therefore $\psi(t)\leq \psi(0)$, which is exactly $\eqref{intertwining}$.\\
We are interested in the behavior of the Poincaré constant along the Polchinski flow. The Poincaré constant $C_P(\nu_t)$ of $\nu_t$ is defined as the smaller constant $K_t$ such that for all $\nu_t$-centered and compactly supported $\phi$, $$\int_{R^d} \phi^2d\nu_t \leq K_t \int_{\R^d} |\nabla \phi|_{C_t'}^2d\nu_t. $$
Contrary to the usual case of study of the Poincaré constant along a flow presented in Section $\ref{sectiongammacalc}$, the metric used in the definition of the Poincaré constant now depends on time. This is because the Polchinski flow gets smaller and smaller until it reaches the Dirac delta function, so in order to be able to see anything at all, we need to expand the metric by an adequate factor. The adequate factor is $C_t'^{-1}$. In the case of the heat kernel decomposition, $C_t'^{-1}= e^{tC_\infty^{-1}}$ and in the case of the Pauli-Villars decomposition, $C_t'^{-1}=t^2\left(A+\frac{1}{t}\right)^{2}$, which both tend towards $+\infty$ when $t$ goes to $+\infty$. 
By starting the Polchinski flow from time $s$ instead of time zero, we can see by imitating the argument for $\nu_0$ that the measure $\nu_s$ satisfies a Poincaré inequality with a sharp constant 
\begin{equation}\label{poincarpournus}
C_P(\nu_s)\leq |C_s'|\int_s^\infty e^{-2\lambda_t}dt.
\end{equation}
In particular, this implies that $$C_P(\nu_s)\underset{s\rightarrow +\infty}{\longrightarrow} C_P(\nu_\infty)=C_P(\delta_0)=0.$$
To better understand the decay of the Poincaré constant along the flow $(\nu_t)_t$, we need to define the following quantity. 
\begin{defi}\label{defalpha}
For all $t\geq 0$, we define $\alpha_t'$ as the biggest eigenvalue of the matrix $$(C_t')^\frac{1}{2} \left(\nabla^2V_t + (C_\infty -C_t)^{-1}\right)(C_t')^\frac{1}{2},$$ where $(C_t')^\frac{1}{2}$ denotes the positive semidefinite square root of the positive semidefinite matrix $C_t'$.
\end{defi}
The quantity $\alpha_t'$ is defined in order to have the following inequality:
$$(C_t')^\frac{1}{2} \left(\nabla^2V_t + (C_\infty -C_t)^{-1}\right)(C_t')^\frac{1}{2}\leq \alpha_t'\, I_d, $$ where $I_d$ is the identity matrix on $\R^d$. Since $C_t$ increases from zero to $C_\infty$, the $\alpha_t'$ will be positive as soon as $\nabla^2 V_t$ is not too negative. In particular, this is so under the Bakry-Emery multiscale criterion when the $\lambda_t'$ are not too often negative. But this is a necessary condition for the Bakry-Emery multiscale criterion to give a non-trivial Poincaré inequality, because of the necessity of the integrability of $e^{-2\lambda_t}$ at $+\infty$ (see Formula $\eqref{piundermbe}$). In the log-concave case, i.e. when $\nabla V_0\geq 0$, then with the heat-kernel covariance decomposition $C_t=C_\infty-C_\infty e^{-tC_\infty^{-1}}$, we get that $\alpha_t'\geq \lambda_{\mathrm{max}}(C_\infty^{-1})>0$.

\subsection{Main result}

We can now state the main theorem of this paper.
\begin{thm}\label{thmmonotonicity}
If the distribution $\nu_0$ given in $\eqref{nuzero}$ is such that its Polchinski flow $(\nu_t)_t$ satisfies the multiscale Bakry-Emery criterion $\eqref{multiscaleeq}$, then for all $0\leq s\leq t$, we have
\begin{equation}\label{ineqthmprincipal}
C_P(\nu_s)\leq e^{(\alpha_t-\alpha_s)-2(\lambda_t - \lambda_s)}C_P(\nu_t), 
\end{equation}
where $\lambda_t := \int_0^t \lambda_s'ds$ is defined from the multiscale curvatures $\lambda_t'$, and $\alpha_t := \int_0^t \alpha_s'ds$ are given in Definition $\ref{defalpha}$. Moreover, if the distribution $\nu_0$ has discrete spectrum (in the sense of improved Poincaré inequalities), then for all $0\leq s\leq t$,
\begin{equation}\label{higherineq}
\lambda_k(\nu_t)\leq e^{(\alpha_t-\alpha_s)-2(\lambda_t - \lambda_s)}\lambda_k(\nu_s), 
\end{equation}
where $\lambda_k$ denotes the $k$-th non-zero eigenvalue.
\end{thm}
Let us underline the following.
\begin{enumerate}
\item When $\nu_0=\gamma$ (i.e. $C_\infty=I_d$ and $V_0=0$), then $\nu_t=\gamma_{e^{-t}I_d}$, so $C_p(\nu_t)=e^{-t}$, $\lambda_t=\alpha_t=t$ and $\eqref{ineqthmprincipal}$ is trivially true.
\item When $\nu_0=\gamma_{C_\infty}$ (i.e. $V_0=0$), then $\nu_t=\gamma_{C_\infty e^{-tC_\infty^{-1}}}$, so $C_p(\nu_t)=\lambda_{\mathrm{max}}(C_\infty e^{-tC_\infty^{-1}})$, $\lambda_t=\frac{\lambda_{\mathrm{min}}(C_\infty)}{2}t$, $\alpha_t=\lambda_{\mathrm{max}}(C_\infty)t $ and $\eqref{ineqthmprincipal}$ becomes $$\lambda_{\mathrm{max}}(C_\infty e^{-sC_\infty^{-1}}) \leq e^{(t-s)\left(\lambda_{\mathrm{max}}(C_\infty)-\lambda_{\mathrm{min}}(C_\infty)\right)}\lambda_{\mathrm{max}}(C_\infty e^{-tC_\infty^{-1}}). $$
\item When $\nu_0$ is log-concave, then $\eqref{ineqthmprincipal}$ is a monotonicity type result of the Poincaré constant along a flow of log-concave measures. Contrary to \cite[Theorem 1.1]{KP} which shows the monotonicity along the heat flow, here it is not an exact monotonicity, and moreover the Polchinski semigroup $(P_{0t})_t$ never coincides with the heat semigroup. However, Theorem $\ref{thmmonotonicity}$ has a broader range of application since it applies to probability distributions which are not log-concave, for instance the $\varphi^4$ measures (see Section $\ref{sectionphi4}$).
\end{enumerate}
The sequel is devoted to the proof of Theorem $\ref{thmmonotonicity}$. This proof is inspired by the corresponding proof of \cite[Theorem 1.1]{KP} by B. Klartag and E. Putterman. We begin with the following lemma.

\begin{lemm}\label{lemmdecrraylegh}
Let $\phi:\R^d\rightarrow\R$ be smooth and compactly supported. We denote by $\phi_t:=P_{0,t}\phi$.
Then the Raylegh quotient $$R_\phi(t) := \frac{\E_{\nu_t} |\nabla \phi_t|^2_{C_t'}}{\E_{\nu_t}\phi_t^2} $$ satisfies for all $t>s$, $$R_\phi(t)\leq R_\phi(s)\exp\left((\alpha_t-\alpha_s)-2(\lambda_t - \lambda_s)\right), $$
where the $\lambda_t$ are given by Assumption $\ref{cdgenralized}$, and $\alpha_t'$ is the largest eigenvalue of the symmetric matrix $C_t'\left(\nabla^2V_t + (C_\infty-C_t)^{-1}\right)$.
\end{lemm} 

\begin{dem}
First, let us underline that $\E_{\nu_t}\phi_t = 0$, justifying that $\E_{\nu_t}\phi_t^2$ is the $\nu_t$-variance of $\phi_t$.
We can compute that (see \cite[Proposition 2.1]{BB})
$$\frac{\partial}{\partial t}\E_{\nu_t} |\nabla \phi_t|^2_{C_t'} = \E_{\nu_t}\left[ 2<\nabla L_t\phi_t,\nabla\phi_t>_{C_t'} - L_t|\nabla \phi_t|^2_{C_t'} + |\nabla \phi_t|^2_{C_t''}\right] $$
and
$$\frac{\partial}{\partial t}\E_{\nu_t}\phi_t^2 = -\E_{\nu_t} |\nabla \phi_t|^2_{C_t'}. $$
Hence 
\begin{align*}
\left(\E_{\nu_t}\phi_t^2\right)^2 R'(t) = \E_{\nu_t}\phi_t^2 \E_{\nu_t}\left[ 2<\nabla L_t\phi_t,\nabla\phi_t>_{C_t'} - L_t|\nabla \phi_t|^2_{C_t'} + |\nabla \phi_t|^2_{C_t''}\right] + \left( \E_{\nu_t} |\nabla \phi_t|^2_{C_t'}\right)^2.
\end{align*}
But the Bochner-type formula for the operator $\Gamma_{2,t}$ associated to the diffusion operator $L_t$ gives 
$$-2\,\Gamma_{2,t}(\phi_t)= 2<\nabla L_t\phi_t,\nabla\phi_t>_{C_t'} - L_t|\nabla \phi_t|^2_{C_t'} =  - ||\nabla^2\phi_t||_{C_t'}^2 - 2<\nabla^2V_t C_t' \nabla\phi_t, \nabla\phi_t>_{C_t'}.$$
So using following multiscale Bakry-Emery condition $\eqref{multiscaleeq}$ (see Assumption $\ref{cdgenralized}$),
$$ C_t'\,\nabla^2 \,V_tC_t' \geq \frac{1}{2}C_t'' + \lambda_t' C_t',$$
one gets $$- 2<\nabla^2V_t C_t' \nabla\phi_t, \nabla\phi_t>_{C_t'} + |\nabla \phi_t|^2_{C_t''} \leq -2\lambda_t' |\nabla \phi_t|_{C_t'}^2.$$
Therefore we obtain
\begin{align*}
R'(t) & \leq \frac{\E_{\nu_t}\left(\phi_t^2\right)E_{\nu_t}\left[- ||\nabla^2\phi_t||_{C_t'}^2-2\lambda_t' |\nabla \phi_t|_{C_t'}^2 \right] + \left( \E_{\nu_t} |\nabla \phi_t|^2_{C_t'}\right)^2}{\left(\E_{\nu_t}\phi_t^2\right)^2} \\
& =   \frac{ \left( \E_{\nu_t} |\nabla \phi_t|^2_{C_t'}\right)^2 -\E_{\nu_t}\left(\phi_t^2\right)E_{\nu_t}\left( ||\nabla^2\phi_t||_{C_t'}^2\right) }{\left(\E_{\nu_t}\phi_t^2\right)^2} - 2\lambda_t' R_\phi(t)
\end{align*}
Moreover, by denoting by $\Gamma_2^{\cL_t} $ the $\Gamma_2$-operator of $\cL_t$,  we have
\begin{align*}
\E_{\nu_t} |\nabla \phi_t|^2_{C_t'} &= -\E_{\nu_t}\left(\phi_t\cL_t\phi_t\right) \\
& \leq \left( \E_{\nu_t}\phi_t^2\right)^\frac{1}{2} \left(  \E _{\nu_t}\left((-\cL\phi_t)^2\right)\right)^\frac{1}{2}\\
& = \left( \E_{\nu_t}\phi_t^2\right)^\frac{1}{2} \left( \E_{\nu_t}\Gamma_2^{\cL_t}(\phi_t) \right)^\frac{1}{2}\\
& = \left( \E_{\nu_t}\phi_t^2\right)^\frac{1}{2} \left( E_{\nu_t}\left( ||\nabla^2\phi_t||_{C_t'}^2\right) + \E_{\nu_t}\left(<\left(\nabla^2V_t +(C_\infty - C_t )^{-1} \right)C_t' \nabla\phi_t, \nabla\phi_t>_{C_t'} \right) \right)^\frac{1}{2}. 
\end{align*}
Since by definition we have $C_t'\left(\nabla^2V_t + (C_\infty -C_t)^{-1}\right)C_t'\leq \alpha_t'C_t'$, we obtain  
\begin{align*}
\left( \E_{\nu_t} |\nabla \phi_t|^2_{C_t'}\right)^2 -\E_{\nu_t}\left(\phi_t^2\right)E_{\nu_t}\left( ||\nabla^2\phi_t||_{C_t'}^2\right) & \leq \E_{\nu_t}\left(\phi_t^2 \right) \E_{\nu_t}\left(<\left(\nabla^2V_t +(C_\infty - C_t )^{-1} \right)C_t' \nabla\phi_t, \nabla\phi_t>_{C_t'} \right)\\
& \leq \alpha_t' \,\E_{\nu_t}\left(\phi_t^2 \right)  \E_{\nu_t}|\nabla \phi_t|^2_{C_t'},
\end{align*}
where $\alpha_t:=\int_0^t \alpha_s'ds$.
Finally, the derivative of $R_\phi$ satisfies
$$R_\phi'(t)\leq \left(\alpha_t' -2\lambda_t'\right)R_\phi(t),$$
and the Gronwall lemma gives
$$R_\phi(t)\leq R_\phi(s)\exp\left((\alpha_t-\alpha_s)-2(\lambda_t - \lambda_s)\right), $$
where $\alpha_t:=\int_0^t \alpha_s'ds$. The proof is complete.
\end{dem}
The Cauchy-Schwarz type inequality $$\forall u\in C^2_c(\R^d),\quad \left( \E_{\nu_t} |\nabla u|^2_{C_t'}\right)^2 \leq \E_{\nu_t}\left(u^2\right)E_{\nu_t}\left( ||\nabla^2 u||_{C_t'}^2\right) $$ never holds true. Indeed, coupled with the Poincaré inequality, it implies
$$\E_{\nu_t}\left(u^2\right)\leq C_P(\nu_t)E_{\nu_t}\left( ||\nabla^2 u||_{C_t'}^2\right), $$ which is easily seen to be false by taking $u$ as an affine function. However, one can see in the proof of Lemma $\ref{lemmdecrraylegh}$ that the Raylegh quotient satisfies the following non linear differential inequation $$R_\phi'(t)\leq R_\phi(t)^2 -\lambda_t'R_\phi(t).$$ But it is not clear that this inequation gives more information about $R_\phi$ than we got in the proof.
We can now prove Theorem $\ref{thmmonotonicity}$.\\
\begin{demof}{\textbf{Theorem $\ref{thmmonotonicity}$}}
By definition of the spectral gap, one has $$\frac{1}{C_P(\nu_0)} = \inf \{ R_\phi(0)\,|\, \phi\in C_c^\infty(\R^d) \}, $$ where $C_c^\infty(\R^d)$ denotes the set of all smooth compactly supported functions on $\R^d$. Let then $t\geq 0$. For all $\varepsilon>0$, there exists $\phi\in C_c^\infty(\R^d)$ such that $$R_\phi(0)<\frac{1}{C_P(\nu_0)}+\varepsilon.$$
By Lemma $\ref{lemmdecrraylegh}$, it follows that
$$e^{2\lambda_t-\alpha_t}R_\phi(t)<\frac{1}{C_P(\nu_0)}+\varepsilon. $$ But by the definition of the Poincaré constant, $\frac{1}{C_P(\nu_t)}\leq R_\phi(t)$, so we get $$\forall \varepsilon>0,\quad \frac{e^{2\lambda_t-\alpha_t}}{C_P(\nu_t)} < \frac{1}{C_P(\nu_0)}+\varepsilon, $$ from which we get 
$$C_P(\nu_t)\geq e^{2\lambda_t-\alpha_t}C_P(\nu_0) $$ by let $\varepsilon$ go to zero. Inequality $\eqref{ineqthmprincipal}$ is therefore proven thanks to the semigroup property: the idea is to start the flow from the time $s$ instead of $0$, and with the same reasoning, we obtain exactly $\eqref{ineqthmprincipal}$.\\
Assume that the distribution $\nu_0$ has discrete spectrum (in the sense of improved Poincaré inequalities), and let us denote by $\lambda_k(\nu_0)$ its $k$-th non-zero eigenvalue. By the min-max theorem, there exists a ($k+1$)-dimensional subspace $E\subset H^1(\nu_0)$ such that $\forall \phi \in E\setminus \{0\}$, $R_\phi(0)\leq \lambda_k(\mu)$. Let us then define the space $E_t:=\{Q_t(\phi)\,|\,\phi\in E\}$. Since for all $t\geq 0$, the operator $P_{0t}(f) = e^{V_t}\gamma_{C_t}\ast(fe^{-V_0})$ is a bijection from $L^1(\nu_0)$ onto $L^1(\nu_0)$, (we can recover the Fourier transform of $fe^{-V_0}$ from the Gaussian convolution $\gamma_{C_t}\ast(fe^{-V_0})$), then $E_t$ is a ($k+1$)-dimensional vector space in $H^1(\mu_t)$. In addition, by using Lemma $\ref{lemmdecrraylegh}$, we have that for all $\psi=P_{0t}(\phi)\in E_t$, $$\frac{\E_{\nu_t} |\nabla \psi|^2_{C_t'}}{\E_{\nu_t}\left[\psi^2\right]} = R_\phi(t) \leq e^{\alpha_t-2\lambda_t}R_\phi(0)\leq e^{\alpha_t-2\lambda_t} \lambda_k(\nu_0).$$ Again by the min-max theorem, we can therefore conclude that $$\lambda_k(\nu_t)\leq e^{\alpha_t-2\lambda_t}\lambda_k(\nu_0),$$ from which we can derive Inequality $\eqref{higherineq}$ by the semigroup property.
\end{demof}
\section{Application to the $\varphi^4$ models}\label{sectionphi4}

In this section, we illustrate Theorem $\ref{thmmonotonicity}$ on the general $\varphi^4$ models (see \cite[Chapter 1]{bauerschmidt2019introduction}) for a more detailed introduction to this model from Quantum Field Theory). Let $\Lambda$ be a finite set, $A=(A_{x,y})_{x,y\in\Lambda}$ be a symmetric positive definite matrix, and let $g>0$, $\nu\in\R$ and $h\in\R^\Lambda$ be constants. The $\varphi^4$ general measure is a probability distribution on $\R^\Lambda$ and is given by the density
\begin{equation}\label{phi4}
d\mu_{A,g,\nu}^{\Lambda,h}(\varphi) = \frac{1}{Z}\exp\left(-\frac{1}{2}(\varphi,A\varphi) -V_0(\varphi) + (h,\varphi) \right)d\varphi
\end{equation} 
where $Z$ is a normalization constant, $(\cdot,\cdot)$ denotes the dot product on $\R^\Lambda$, and $$V_0(\varphi)=\sum_{\Lambda}\left( \frac{1}{4}g\varphi_x^4 + \frac{1}{2}\nu\varphi_x^2\right). $$
The probability distribution $\mu_{A,g,\nu}^{\Lambda,h}$ is of the form $\eqref{nuzero}$ with $C_\infty =A^{-1}$. We use the Pauli-Villars covariance regularization $C_0=0$ and $C_t=\left(A+\frac{1}{t}\right)^{-1}$, $t>0$. Then the multiscale Bakry-Emery criterion $\eqref{multiscaleeq}$ is satisfied with 
\begin{equation}\label{lambdaphi4}
\lambda_t'=\frac{1}{t}-\frac{\chi_t}{t^2},\quad t>0,
\end{equation} where the quantity $\chi_t$ is the susceptibility, and is defined as the maximum over $\Lambda$ of the correlations between two points:
$$\chi_t:= \max_{x\in\Lambda}\sum_{y\in\Lambda}\int_{\R^\Lambda}\varphi_x\varphi_y\, d\mu_{A,g,\nu+\frac{1}{t}}^{\Lambda,0}(\varphi). $$ 
This was proved in \cite[Proposition 2.5]{bauerschmidt2022log} and used to show that the $\varphi^4$ measures satisfy a logarithmic-Sobolev inequality. Therefore, Theorem $\ref{thmmonotonicity}$ applies to the measure $\mu_{A,g,\nu}^{\Lambda,h}$, and we obtain the following.
\begin{prop}
Let $\nu_t$ be the Polchinski flow associated to the distribution $\mu_{A,g,\nu}^{\Lambda,h}$. Then the Poincaré constant along this flow satisfies 
\begin{equation}\label{poincaremodelphi4}
\forall 0\leq s<t,\quad C_P(\nu_s)\leq e^{(\alpha_t-\alpha_s)-2(\lambda_t - \lambda_s)}C_P(\nu_t),
\end{equation} with $\lambda_t:=\int_0^t\left(\frac{1}{s}-\frac{\chi_s}{s^2}\right)ds$, and $\alpha_t:=\int_0^t\alpha_s'ds$ where
\begin{equation}\label{alphatphi4}
\alpha_t'=\frac{1}{t} -\frac{1}{t^2}\underset{{\varphi\in\R^\Lambda}}{\inf} \lambda_{\min}(\Sigma_t(\varphi))+\lambda_{\max}(A)(t\lambda_{\max}(A)+1),
\end{equation} where $\Sigma_t(\varphi)$ denotes the covariance matrix of the probability distribution $\mu_{A,g,\nu+\frac{1}{t}}^{\Lambda,C_t^{-1}\varphi}$.
\end{prop} 
\begin{dem}
It only remains to prove the bound for $\alpha_t'$. On the one hand, we compute
$(C_\infty -C_t)^{-1} = A(tA+I_d) $ and moreover $(C_t')^\frac{1}{2}=\frac{C_t}{t}$ decreases in t, so $(C_t')^\frac{1}{2}\leq (C_0')^\frac{1}{2}=I_d$. On the other hand, by using the same calculation as in \cite[Section 2]{bauerschmidt2022log}, we have $\nabla^2V_t(\varphi) = C_t^{-1}-C_t^{-1}\Sigma_t(\varphi) C_t^{-1} $, where $\Sigma_t(\varphi)$ denotes the covariance matrix of the probability distribution $\mu_{A,g,\nu+\frac{1}{t}}^{\Lambda,C_t^{-1}\varphi}$.
Therefore we obtain that for all $\varphi\in\R^\Lambda$,
\begin{align*}
(C_t')^\frac{1}{2} \left(\nabla^2V_t(\varphi) + (C_\infty -C_t)^{-1}\right)(C_t')^\frac{1}{2} &= (C_t')^\frac{1}{2} \left( C_t^{-1}-C_t^{-1}\Sigma_t(\varphi) C_t^{-1} + (C_\infty -C_t)^{-1}\right)(C_t')^\frac{1}{2} \\
& \leq \frac{1}{t} -\frac{\lambda_{\min}(\Sigma_t(\varphi))}{t^2}+\lambda_{\max}(A)(t\lambda_{\max}(A)+1),
\end{align*}
where inequations are taken in the sense of quadratic forms.
\end{dem}

In Formula $\eqref{alphatphi4}$, the term $\frac{1}{t} -\frac{1}{t^2}\underset{{\varphi\in\R^\Lambda}}{\inf} \lambda_{\min}(\Sigma_t(\varphi))$ controls the variation of $C_P(\nu_t)$ when $t$ is small, and the term $\lambda_{\max}(A)(t\lambda_{\max}(A)+1)$ controls the variation of $C_P(\nu_t)$ when $t$ is large. If we take $s=0$ and if we let $t$ go to infinity in $\eqref{poincaremodelphi4}$, according to $\eqref{poincarpournus}$, we have that $C_P(\nu_t)$ tends towards zero. However, the exponential term then behaves as its dominant term $e^{t\lambda_{max}(A)^2}$ and thus counterbalances the decrease of $C_P(\nu_t)$.\\

\textbf{Acknowledgments:} {This work was supported by the Labex Cimi of the University of Toulouse. I would like to thank Yair Shenfeld for our interesting conversation which enlightened me for Section $\ref{sectiontransport}$ and for his comments on the first version of this paper. }

\bibliographystyle{plain}
\bibliography{mabibliographie}
\end{document}